\documentclass[a4paper, 10pt]{article}

\usepackage[numbers]{natbib}

\usepackage{filecontents}
\usepackage{float}

\usepackage[english]{babel}
\usepackage[utf8]{inputenc}
\usepackage[T1]{fontenc}

\usepackage{amsmath}
\usepackage{amssymb}

\usepackage{xcolor}
\usepackage{bbm}

\usepackage{enumerate}

\usepackage{scalerel} 
\usepackage{subfigure}

\usepackage{graphicx}

\usepackage{amsthm}
\newtheorem{theorem}{Theorem}
\numberwithin{theorem}{section}
\newtheorem{lemma}[theorem]{Lemma}
\newtheorem{prop}[theorem]{Proposition}

\newtheorem{corollary}[theorem]{Corollary}

\theoremstyle{definition}

\newtheorem{definition}[theorem]{Definition}

\newtheorem{remark}[theorem]{Remark}

\usepackage{pgfplots}

\pgfplotsset{compat = newest}

\newcommand{\N}{\mathbb{N}}
\newcommand{\Q}{\mathbb{Q}}
\newcommand{\R}{\mathbb{R}}

\newcommand{\E}{\mathbb{E}}

\newcommand{\A}{\mathcal{A}}

\newcommand{\F}{\mathcal{F}}
\newcommand{\G}{\mathcal{G}}
\newcommand{\B}{\mathcal{B}}

\newcommand{\PP}{\mathcal{P}}

\newcommand{\ind}[1]{{\textbf{1}}_{{#1}}}

\newcommand{\id}{\textup{id}}

\newcommand{\cpl}{\textup{Cpl}}

\newcommand{\W}{\mathcal{W}}

\renewcommand{\epsilon}{\varepsilon}
\renewcommand{\subset}{\subseteq}

\renewcommand{\P}{\mathbb{P}}

\newcommand{\proj}{\operatorname{proj}}

\let\oldmarginpar\marginpar
\renewcommand\marginpar[1]{\-\oldmarginpar[\raggedleft\footnotesize #1]{\raggedright\footnotesize\color{red} #1}}
\usepackage[hidelinks]{hyperref}

\usepackage{tikz}
\usepackage{subfigure}
\usepackage{pgfplots}
\pgfplotsset{compat=newest}

\title{Measurable Universality and Selection in Optimal Transport}

\author{Daniel Toneian\footnote{University of Vienna, Faculty of Mathematics: daniel.toneian@univie.ac.at}}

\begin{document}
	\maketitle
	\begin{abstract}
    We study the dependence of optimal transport problems on their marginals.
    We show that parametrized Borel minimization over an uncountable Polish
    space can be represented by a fixed Borel transport cost and continuously
    varying equal marginals. The Kantorovich value is universally measurable
    for Borel costs, but on every uncountable Polish space it can fail to be
    Borel measurable even for costs taking only the values zero and one and
    equal marginals.

    For lower semicontinuous costs, we show how to select
    optimal plans measurably and, when Kantorovich-optimal maps exist, how to
    choose them jointly measurably. We also obtain measurable selection
    results for Wasserstein barycenters and optimal Skorokhod embeddings.

	\end{abstract}

\section{Introduction}

We study how the value and optimizers of an optimal transport problem
depend on its marginals. Our main result is that optimal transport with
Borel costs can represent arbitrary parametrized Borel minimization
problems. A particular consequence is that for Borel
costs the optimal transport value can fail to be Borel measurable,
even when the cost takes only the values zero and one and the two
marginals are equal.

This dependence on the marginals matters, for instance, when transport is conditioned
on observed information. The two conditional laws are probability
kernels indexed by the observation. To solve the transport problem
conditionally, we need the value to depend measurably on the
observation and we need to choose optimal plans measurably.

Let $X$ and $Y$ be Polish spaces. For a Polish space $Z$,
write $\PP(Z)$ for the space of Borel probability measures on $Z$,
equipped with the topology of weak convergence. For $\mu\in\PP(X)$ and
$\nu\in\PP(Y)$, write $\cpl(\mu,\nu)$ for the couplings, or transport
plans, of $\mu$ and $\nu$. These are the laws $\pi\in\PP(X\times Y)$
whose marginals satisfy $\pi_X=\mu$ and $\pi_Y=\nu$. For a Borel cost
$c:X\times Y\to[0,\infty]$, the Kantorovich value is
\begin{equation}\label{intro:transport}
    V_c(\mu,\nu):=\inf_{\pi\in\cpl(\mu,\nu)}\int_{X\times Y}c\,d\pi.
\end{equation}
A coupling attaining this infimum is called an optimal transport plan.

We call the representation result \emph{measurable universality}.

\begin{theorem}[Measurable universality]\label{intro:universality}
For every uncountable Polish space $M$ there exist $\mu\in\PP(M)$ and
a Borel map $D:M\times M\to M$ such that, for every Polish space $P$
and every Borel $g:P\times M\to[0,\infty]$,
\[
    G(p):=\inf_{z\in M}g(p,z)
    =\inf_{\pi\in\cpl(\mu,\mu)}
       \int_{M\times M} g\bigl(p,D(x,y)\bigr)\,\pi(dx,dy),\qquad p\in P.
\]
There is also a single Borel cost $\tilde c$ on $(P\times M)^2$ such
that $G(p)=V_{\tilde c}(\eta_p,\eta_p)$ for every $p\in P$, where
$\eta_p=\delta_p\otimes\mu$ depends continuously on $p$. If $g$ is
bounded, $\tilde c$ can be chosen bounded.
\end{theorem}

The idea is to couple a uniform point on the circle with its translate
by $z$. Both marginals remain uniform, and their difference modulo one
is $z$, so every competitor is realized by a coupling. A Borel
isomorphism transfers this construction to $M$. The
two representations are proved in Theorem~\ref{prop:universality}
and Corollary~\ref{cor:universality}.

\begin{theorem}[Measurability of transport values]\label{intro:nonborel}
For every Borel cost $c:X\times Y\to[0,\infty]$ on Polish spaces,
$V_c$ is universally measurable. For every uncountable Polish space
$X$, there is a Borel cost $c:X\times X\to\{0,1\}$ such that
$\mu\mapsto V_c(\mu,\mu)$ is not Borel measurable.
\end{theorem}

Universal measurability here means that for any Borel
probability law on the space of pairs of marginals, the value is
measurable once all subsets of Borel sets of probability zero are
included among the measurable sets. The proof of the first claim of Theorem~\ref{intro:nonborel} is in
Lemma~\ref{lemma:universal_value}, and the descriptive set theory used
there is reviewed in the Appendix. For the second claim, take $g$ from Theorem~\ref{intro:universality} to be zero on a Borel
set with non-Borel projection and one elsewhere. This yields the binary cost
example in Theorem~\ref{theo:OT_not_measurable}.
Even a closed set can have a non-Borel projection \citep{Ke95}, so a
lower semicontinuous objective can yield a non-Borel $G$. Such a $G$
cannot arise from a lower semicontinuous transport cost with
continuously varying marginals.

We next turn to the selection of optimizers.

\begin{theorem}[Selection of optimal plans and maps]\label{intro:selection}
Let $c:X\times Y\to[0,\infty]$ be lower semicontinuous, with $X$ and
$Y$ Polish. Then $V_c$ is lower semicontinuous and the defining infimum is
attained for every pair of marginals. Let $(E,\mathcal E,m)$ be a
$\sigma$-finite measure space and let $e\mapsto(\mu_e,\nu_e)$ be an
$\mathcal E$-measurable family of marginals. There is an
$\mathcal E$-measurable choice $e\mapsto\pi_e$ of plans that are
optimal for $m$-almost every $e$.
If each pair $(\mu_e,\nu_e)$ admits a Kantorovich-optimal plan induced
by a map, there is an $\mathcal E\otimes\B(X)$-measurable map
$T:E\times X\to Y$ such that $\bigl(\mathrm{id}_X,T(e,\cdot)\bigr)_\#\mu_e$ is
an optimal transport plan with marginals $(\mu_e,\nu_e)$ for
$m$-almost every $e$.
\end{theorem}

Attainment and lower semicontinuity follow from
compactness and lower semicontinuity of the cost integral
\citep{Vi09}. See Theorem~\ref{theo:selection_kantorovich} and the result
on maps that follows it for the selection statements.
The selection argument applies the
Jankov--von Neumann theorem to the set of admissible minimizers.
We also give a further result for Borel costs. Measurable selection is possible provided
optimal plans exist for $m$-almost every parameter
(Theorem~\ref{theo:selection_determinacy}).

The same selection principle gives measurable choices of Wasserstein
barycenters for measures supported on a common compact convex set
(Theorem~\ref{theo:selection_barycenters}). In the transport formulation
of Beiglb\"ock, Cox, and Huesmann \citep{BeCoHu17}, we also select
optimal randomized Skorokhod embeddings measurably for centered target
laws with finite second moment and lower semicontinuous costs. With
an independent uniform random variable available at time zero, we
obtain a measurable choice of ordinary optimal stopping times
(Theorems~\ref{theo:skorokhod} and~\ref{theo:selection_stopping_times}).
\par

Bogachev and Malofeev \citep[Theorem~4.2 and Corollary~4.4]{BoMa20}
prove Borel selection of optimal plans and jointly Borel
disintegration for Souslin parameter spaces and finite optimal values,
allowing both marginals and lower semicontinuous costs to vary with
the parameter. Leskelä and Vihola \citep{LeVi17} prove Borel selection of optimal plans on $\R^d$ for lower semicontinuous costs satisfying a linear lower bound and measures with finite first moments. Fontbona, Guérin and Méléard \citep{FoGuMe10} establish joint measurability of optimal maps on $\R^d$ assuming that each transport problem admits a unique optimal transport plan and that this plan is induced by a map. Several of the results in this paper have appeared in a preliminary form in the author's thesis \citep[Chapters~3--6]{To19}.
\par

\section{Measurable Universality of Optimal Transport}

We consider parametrized optimization problems of the form
$$G(p):= \inf_{x \in M} g(p,x),$$
where $P$ and $M$ are Polish spaces and the function $g: P \times M \to [0,+\infty]$ is Borel measurable. The next theorem will show that, as long as $M$ is uncountable, any such problem has a representation as an optimal transport problem.

\begin{theorem}\label{prop:universality}
    Let $M$ be an uncountable Polish space. Then there exists a Borel function 
    $$D: M \times M \to M$$
    and a probability measure $\mu \in \PP(M)$ so that for every Polish space $P$ and every Borel function
    $$g: P \times M \to [0,+\infty]$$
    one has, for every $p \in P$,
    $$G(p):=\inf_{x \in M} g(p,x) = \inf_{\pi \in \cpl(\mu,\mu)}\int_{M \times M} g\bigl(p,D(x,y)\bigr) \, \pi(dx,dy).$$
\end{theorem}

\begin{proof}
    For any $x\in \R$ we denote the fractional part $\{x\}$ of $x$ by $\{x\}:=x - \lfloor x \rfloor$, where $\lfloor x\rfloor$ is the floor of $x$. We separate the proof into two steps.

    For the first step, assume $M=[0,1)$.   Let $\mu:=\lambda$ be the Lebesgue measure on $[0,1)$ and $D(x,y)=\{y-x\}$.
    Let $p \in P$ and $z\in [0,1)$ be given. We show that there exists a coupling $\pi^z \in \cpl(\mu,\mu)$ so that
    $$ g(p,z)=\int g\bigl(p,D(x,y)\bigr) \, \pi^z(dx,dy).$$
    We first introduce the map $\tau_z$ as
    $$ \tau_z:  [0,1) \to [0,1) \qquad \tau_z(x):= \{x+z\}$$
    and then define the probability measure $\pi^z\in \PP([0,1) \times [0,1))$ as
    $$ \pi^z:= (\id,\tau_z)_{\#}\mu,$$
    which puts all its mass on the diagonal of $[0,1)^2$ shifted by $z$. Note that
    $$ (\tau_z)_{\#}\mu=\mu,$$
    because translation modulo $1$ preserves the Lebesgue measure, hence $\pi^z \in \cpl(\mu,\mu).$  Since $\bigl\{\{x+z\}-x\bigr\}=z$ whenever $0\leq z <1$ and $x\in \R$, we see that
    \begin{align*}  
    &\int_{M\times M} g\bigl(p,D(x,y)\bigr) \, \pi^z(dx,dy) = \int_{[0,1)\times [0,1)} g\bigl(p,\{y-x\}\bigr) \, \pi^z(dx,dy) \\
    = &\int_{[0,1)} g\Bigl(p,\bigl\{\{x+z\}-x\bigr\}\Bigr)\, \lambda(dx) = \int_{[0,1)} g(p,z) \, \mu(dx)= g(p,z).
        \end{align*}
    As $z$ was arbitrary, we get
    \begin{equation*}
         \inf_{\pi \in \cpl(\mu,\mu)}\int_{M \times M} g\bigl(p,D(x,y)\bigr) \, \pi(dx,dy) \leq \inf_{z\in M} g(p,z).
    \end{equation*}
    For the reverse inequality, note that for any $p \in P$ and  $x,y \in M$ we have
    $$\inf_{z \in M } g(p,z) \leq g\bigl(p,D(x,y)\bigr),$$
    and hence
    $$ \inf_{z \in M } g(p,z) \leq \int_{M \times M} g\bigl(p,D(x,y)\bigr) \, \pi(dx,dy)$$
    for any probability measure $\pi \in \PP(M \times M)$. Taking the infimum over $\pi \in \cpl(\mu,\mu)$ gives
    $$  \inf_{z \in M } g(p,z) \leq \inf_{\pi \in \cpl(\mu,\mu)}\int_{M \times M} g\bigl(p,D(x,y)\bigr) \, \pi(dx,dy).$$
    
    For the second step, let $M$ be an arbitrary uncountable Polish space. By the Borel Isomorphism Theorem \cite[Theorem 15.6]{Ke95}, there exists a Borel isomorphism 
    $$\phi: M \to [0,1).$$
    Set
    $$\mu := (\phi^{-1})_{\#}\lambda, $$
    where $\lambda$ is the Lebesgue measure on $[0,1)$, and
    $$D: M \times M \to M, \qquad D(x,y):= \phi^{-1}\bigl(\{\phi(y)-\phi(x)\}\bigr).$$
    Since $\phi, \phi^{-1}$ and the fractional-part map are Borel measurable, so is $D$. Note that once the Borel isomorphism $\phi$ is fixed, both $\mu$ and $D$ are independent of $P$ and $g$.
    Define 
    $$ \tilde g: P \times [0,1) \to [0,+\infty] \qquad \tilde g(p,t):=g(p, \phi^{-1}(t)).$$
    Then $\tilde g$ is Borel measurable. As $\phi$ is bijective,
    $$ \inf_{x\in M}g(p,x)= \inf_{t \in [0,1)} \tilde g(p,t).$$
    Fix $p \in P$. Applying the first step to $\tilde g$, we obtain
    $$\inf_{t\in [0,1)} \tilde g(p,t) = \inf_{\tilde \pi \in \cpl(\lambda, \lambda)} \int \tilde g\bigl(p,\{\tilde y-\tilde x\}\bigr) \, \tilde \pi (d\tilde x, d \tilde y). $$
    Moreover,
    $$\pi \mapsto (\phi \times \phi)_{\#}\pi$$
    is a bijection between $\cpl(\mu,\mu)$ and $\cpl(\lambda,\lambda)$. Therefore
    \begin{align*}
        &\inf_{z \in M} g(p,z)= \inf_{ t \in [0,1)} \tilde g(p, t) = \inf_{\tilde \pi  \in \cpl(\lambda, \lambda)} \int \tilde g\bigl(p,\{\tilde y-\tilde x\}\bigr) \, \tilde \pi (d\tilde x, d \tilde y) \\
        = &\inf_{\tilde \pi \in \cpl(\lambda, \lambda)}\int_{[0,1)\times [0,1)} g\Bigl(p,\phi^{-1}\bigl(\{\tilde y-\tilde x\}\bigr)\Bigr) \, \tilde\pi (d\tilde x,d\tilde y)\\
        = &\inf_{\pi \in \cpl(\mu, \mu)} \int_{M \times M}g\Bigl(p,\phi^{-1}\bigl(\{\phi(y)-\phi(x)\}\bigr)\Bigr) \, \pi (dx,dy) \\
        =&\inf_{\pi \in \cpl(\mu, \mu)} \int_{M \times M}g\bigl(p,D(x,y)\bigr) \, \pi (dx,dy),
    \end{align*}
    which completes the proof.
\end{proof}

In the preceding theorem, $p$ enters the transport cost. The following
corollary encodes $p$ in the marginals and uses a single cost. In
particular, any parametrized minimization problem with a nonnegative
Borel objective and a Borel constraint set can be represented as an
optimal transport problem by assigning cost $+\infty$ to inadmissible
competitors. A representation with continuously varying marginals
cannot always use a lower semicontinuous transport cost, as shown in
Remark~\ref{rem:lsc_universality}.
\par

\begin{corollary}\label{cor:universality}
   Let $P$ and $M$ be Polish spaces, with $M$ uncountable and let the map
   $$g:P \times M \to [0,+\infty]$$
   be Borel measurable. Define
   $$G: P \to [0,+\infty], \qquad G(p):=\inf_{x \in M}g(p,x).$$
   Then there exist a probability measure $\mu \in \PP(M)$, a Borel cost function 
   $$c: (P \times M) ^2 \to  [0,+\infty] $$ 
   and a continuous map
   $$\Phi_\mu: P \to \PP(P \times M) \times \PP(P \times M), \qquad p \mapsto (\delta_p \otimes \mu, \delta_p \otimes \mu),$$
   so that 
   $$G= V_c \circ \Phi_\mu.$$
   Moreover, if $g$ is bounded then $c$ can be chosen bounded.
\end{corollary}

\begin{proof}
    Let $\mu \in \PP(M)$ and $D: M \times M \to M $ be as in Theorem~\ref{prop:universality} and define
    $$ c\bigl((p_1,x),(p_2,y)\bigr) := g(p_1,D(x,y)).$$
    Then $c$ is Borel measurable (and bounded if $g$ is bounded). Fix $p\in P$. As every coupling between $\delta_p \otimes \mu$ and itself is concentrated on 
    $$ (\{p\} \times M) \times (\{p\} \times M), $$
    such couplings can be naturally identified with elements of $\cpl(\mu,  \mu).$ Hence, by Theorem~\ref{prop:universality}, we get
    $$ V_c\bigl(\Phi_\mu(p)\bigr) = \inf_{\pi \in \cpl (\mu, \mu) } \int_{M \times M} g\bigl(p,D(x,y)\bigr) \, \pi(dx,dy) = \inf_{z \in M}g(p,z) = G(p). $$
    It remains to show that $\Phi_\mu $ is indeed continuous. Let $(p_n)_{n \in \N}$ be a sequence in $P$ such that $\lim_{n \to \infty} p_n=p$ for some $p \in P$, and let $f: P \times M \to \R$ be bounded and continuous. Then
    \begin{align*}
       &  \lim_{n\to \infty}\int_{P \times M} f \, d(\delta_{p_n} \otimes \mu) = \lim_{n\to \infty} \int_M f(p_n,x) \, \mu(dx)  \\
       &= \int_M f(p,x) \, \mu(dx) = \int_{P \times M} f \, d(\delta_p \otimes \mu),
    \end{align*}
    by dominated convergence. Thus,
    $$\lim_{n \to \infty} \delta_{p_n}\otimes \mu= \delta_p \otimes \mu$$
    weakly and hence
    $$\lim_{n\to \infty} \Phi_\mu(p_n)= \Phi_\mu(p)$$
    in $\PP(P \times M)^2$.
\end{proof}

For lower semicontinuous costs, the Kantorovich value is lower
semicontinuous in the marginals. For merely Borel costs the value
remains universally measurable (Lemma~\ref{lemma:universal_value}),
but Borel measurability can fail even for binary costs and equal
marginals.
\par

\begin{theorem}\label{theo:OT_not_measurable}
    
    Let $X$ be an uncountable Polish space. There exists a Borel cost
    $c:X\times X\to\{0,1\}$ for which $\mu\mapsto V_c(\mu,\mu)$ on
    $\PP(X)$ is not Borel measurable. In particular, $V_c$ is not
    Borel measurable on $\PP(X)^2$.
    \par
\end{theorem}

\begin{proof}

    Fix uncountable Polish spaces $P$ and $M$ and set $\Omega=P\times M$.
    Choose an analytic non-Borel set $S\subseteq P$. By
    Theorem~\ref{app:equivalent_formulations_of_analytic}, there is a
    Borel set $A\subseteq P\times M$ with $\proj_P(A)=S$. For
    $g=1-\ind{A}$, we have
    \[
        G(p):=\inf_{z\in M}g(p,z)=1-\ind{S}(p),
    \]
    so $G$ is not Borel measurable.

    Use the cost from the proof of Corollary~\ref{cor:universality},
    $\tilde c\bigl((p,x),(q,y)\bigr)=g\bigl(p,D(x,y)\bigr)$. It takes values in $\{0,1\}$.
    The corollary gives a probability measure $\mu\in\PP(M)$ such that
    $G(p)=V_{\tilde c}(\eta_p,\eta_p)$ for $\eta_p=\delta_p\otimes\mu$.
    Since $p\mapsto\eta_p$ is continuous, the diagonal value
    $\eta\mapsto V_{\tilde c}(\eta,\eta)$ cannot be Borel measurable.

    Now let $X$ be any uncountable Polish space. Choose a Borel
    isomorphism $\phi:X\to\Omega$ \citep[Theorem 15.6]{Ke95} and set
    $c=\tilde c\circ(\phi\times\phi)$. This is a Borel cost taking
    values in $\{0,1\}$. For each $\eta\in\PP(\Omega)$, write
    $\nu=(\phi^{-1})_\#\eta$. Pushing forward by $\phi\times\phi$
    gives a bijection between $\cpl(\nu,\nu)$ and $\cpl(\eta,\eta)$
    that preserves costs, so
    \[
        V_c\big((\phi^{-1})_\#\eta,(\phi^{-1})_\#\eta\big)
        =V_{\tilde c}(\eta,\eta).
    \]
    The pushforward $\eta\mapsto(\phi^{-1})_\#\eta$ is Borel because
    its evaluation on each Borel set $B\subseteq X$ is
    $\eta\mapsto\eta(\phi(B))$, and evaluations generate the Borel
    $\sigma$-algebra of $\PP(X)$. Thus Borel measurability of
    $\nu\mapsto V_c(\nu,\nu)$ would imply Borel measurability of
    $\eta\mapsto V_{\tilde c}(\eta,\eta)$. Borel measurability of
    $V_c$ itself would imply Borel measurability on the diagonal,
    since $\nu\mapsto(\nu,\nu)$ is continuous.
\end{proof}

\begin{remark}\label{rem:lsc_universality}
    Theorem~\ref{prop:universality} and Corollary~\ref{cor:universality} show that it is possible to rewrite generic optimization problems in this specific way as optimal transport problems. There cannot, however, be a version of Corollary~\ref{cor:universality} that starts with a lower semicontinuous function $g$ to be optimized and yields a lower semicontinuous cost function $c$ while retaining a continuous parametrization of the marginals. To see this, let $A \subseteq \N^\N \times \N^\N$ be a closed set whose projection onto the first coordinate is not Borel measurable (such sets exist, as can be seen for example in the proof of \cite[14.2]{Ke95}) and $g(p,x)=1-\ind{A}(p,x)$. Since $A$ is closed, $g$ is lower semicontinuous and the optimization problem $G(p)=\inf_{x \in \N^\N}g(p,x)$ is not Borel measurable. If there were a version of Corollary~\ref{cor:universality} giving a lower semicontinuous cost function $c$, then the associated value function would also be lower semicontinuous by Lemma~\ref{lemma:V_is_lsc}, hence Borel measurable, which is a contradiction. 
\end{remark}

\section{Measurable Selection in Optimal Transport}
We use the notation of the introduction for the
marginals $\mu$ and $\nu$, the cost $c$ and the Kantorovich value
$V_c$ in \eqref{intro:transport}.

A measurable map $T : X \to Y$ is said to transport $\mu$ to $\nu$ if
$T_\#\mu=\nu$. The (Monge) optimal transport problem consists in determining
$$\inf_{\substack{T : X \to Y\\ T_\#\mu=\nu}} \int_X c\bigl(x,T(x)\bigr)\,d\mu(x).$$
Any minimizer is called an \emph{optimal transport map}.

Every transport map $T$ induces a transport plan
$(\mathrm{id},T)_\#\mu \in \cpl(\mu,\nu)$.
Consequently, whenever the Monge problem admits a minimizer, the induced
transport plan is admissible for the Kantorovich problem.

We begin by showing that the optimal transport problem, as a function of the marginals, is lower semicontinuous.

\begin{lemma}\label{lemma:V_is_lsc}
    Let $X,Y$ be Polish spaces and $c: X\times Y \to [0,\infty]$ a lower semicontinuous cost function. Then the value function $V_c$ is lower semicontinuous.
\end{lemma}

\begin{proof}
    It is well known that $\pi \mapsto \int c \, d\pi$ is lower semicontinuous (see, for example, \cite[Lemma 4.3]{Vi09}). Now let $(\mu_n,\nu_n)_{n\in \N}$ be a sequence in $\PP(X)\times \PP(Y)$ converging weakly to a pair $(\mu,\nu) \in \PP(X)\times \PP(Y)$. If 
    $$\liminf_n V_c(\mu_n, \nu_n)= +\infty, $$
    there is nothing to prove. Otherwise, choose a subsequence  $(\mu_{n_k}, \nu_{n_k})_{k\in \N}$ of $(\mu_n, \nu_n)_{n \in \N}$  such that
    $$\lim_{k \to \infty} V_c(\mu_{n_k}, \nu_{n_k}) = \liminf_n V_c(\mu_n, \nu_n). $$
    Let $\pi_k \in \cpl(\mu_{n_k}, \nu_{n_k} )$ for $k \in \N$ such that
   $$ \int c \, d\pi_k < V_c(\mu_{n_k}, \nu_{n_k}) + \frac{ 1}{k}.$$
    The measures $ (\pi_k)_{k\in \N}$ are tight (\cite[Lemma 4.4]{Vi09}). Hence, after passing to a further subsequence, $\lim_{k \to \infty} \pi_k =\hat \pi$ weakly  for some $\hat\pi \in \cpl(\mu,\nu)$. By lower semicontinuity

\begin{align*}
V_c(\mu,\nu)&= \inf_{\pi \in \cpl(\mu,\nu)} \int c \, d\pi \leq \int c \, d\hat\pi \leq \liminf_k \int c \, d\pi_k \\
&\leq \liminf_k \left( V_c(\mu_{n_k}, \nu_{n_k}) + \frac{1}{k} \right)= \liminf_n V_c(\mu_n, \nu_n) ,
\end{align*} 
which completes the proof.
\end{proof}

We introduce a short lemma that we will use several times.

\begin{lemma}\label{lemma:couplings_closed}
    Let $X,Y$ be Polish spaces. Then
      $$G= \{ ( \mu, \nu, \pi) \in  \PP(X) \times \PP(Y) \times \PP(X \times Y) : \pi \in \cpl(\mu, \nu)\} $$
      is closed in $\PP(X) \times \PP(Y) \times \PP(X \times Y).$
\end{lemma}

\begin{proof}
    Let $(\mu_n, \nu_n , \pi_n)_{n\in \N}$ be a sequence in $G$ converging weakly in each component to $(\mu, \nu, \pi)$. As the marginal maps
    $$\pi \mapsto \pi_X, \qquad \pi \mapsto \pi_Y$$
    are continuous for the weak topology, we find that 
    $$\pi_X= \lim_{n \in \N} {(\pi_n)_X} = \lim_{n \in \N} \mu_n =\mu $$
     and
     $$\pi_Y= \lim_{n \in \N} {(\pi_n)_Y} = \lim_{n \in \N} \nu_n =\nu. $$
     Thus, 
    \[(\mu, \nu, \pi) \in G. \qedhere\] 
\end{proof}

If the cost function $c$ is not lower semicontinuous but merely Borel, the value function is at least universally measurable. 

\begin{lemma}\label{lemma:universal_value}
     Let $X,Y$ be Polish spaces and $c: X\times Y \to [0,\infty]$ a Borel measurable cost function. Then the value function $V_c$ is universally measurable.
\end{lemma}

\begin{proof}
    We start by noting that 
    $$I_c: \PP(X \times Y) \to [0,+\infty], \qquad \pi \mapsto \int c\, d\pi$$ is Borel. Indeed, if $c = \ind{A}$ for some Borel $A\subseteq X\times Y$ this follows by the Portmanteau theorem and a monotone-class argument. This implies that $I_c$ is Borel for nonnegative Borel simple cost functions $c$. For arbitrary Borel nonnegative $c$, approximate $c$ by a family of simple Borel functions $(c_n)_{n \in\N}$ so that $c_n \nearrow c$. Then $I_c$ is Borel, as by monotone convergence we know that
    $$I_c(\pi)= \int c \, d\pi = \lim_{n \to \infty} \int c_n \, d\pi, $$
    and thus $I_c$ is the limit of Borel functions and hence Borel itself.

    Next, note that by Lemma~\ref{lemma:couplings_closed} the set
    $$G:= \{ ( \mu, \nu, \pi) \in  \PP(X) \times \PP(Y) \times \PP(X \times Y) : \pi \in \cpl(\mu, \nu)\}  $$
    is closed, and thus Borel. 
    Now, fix $a \in \R$. We find that
    $$ \{(\mu, \nu): V_c(\mu, \nu) < a \} = \proj_{\PP(X) \times \PP(Y)} \Bigl(G\cap\bigl\{(\mu,\nu,\pi):I_c(\pi)<a\bigr\}\Bigr).$$
    As both $G$ and $\{ ( \mu, \nu, \pi) :I_c(\pi) < a\}$ are Borel, we conclude that  $\{(\mu, \nu): V_c(\mu, \nu) < a \}$ is the projection of a Borel set, thus analytic. As this is true for any $a \in \R$, the function $V_c$ is lower semianalytic, and in particular universally measurable.
\end{proof}

Next, we introduce the extended cost functional.

\begin{definition}
    Let  $X,Y$ be Polish spaces and let $c: X\times Y \to [0,\infty]$ be a lower semicontinuous cost function. Then the extended cost functional $C: \PP(X) \times \PP(Y) \times \PP(X \times Y) \to [0, +\infty]$ is given by
    \begin{equation*}\label{transport:C}
   C(\mu, \nu, \pi)=
\begin{cases} 
\int_{X\times Y} c(x,y) \, d\pi, &\pi \in \cpl(\mu,\nu), \\
+\infty, &\text{otherwise.}
\end{cases} 
\end{equation*}
\end{definition}
We proceed by showing that the extended cost functional is lower semicontinuous.
\begin{lemma}\label{lemma:extended_cost_is_lsc}
    Let  $X,Y$ be Polish spaces and $c: X\times Y \to [0,\infty]$ be a lower semicontinuous cost function. Then the extended cost functional $C$ is lower semicontinuous.
\end{lemma}
\begin{proof}
    Define
    $$ G:= \{ (\mu, \nu, \pi) \in \PP(X)\times \PP(Y) \times \PP(X\times Y) : \pi \in \cpl(\mu, \nu)\}.$$
    By Lemma~\ref{lemma:couplings_closed} $G$ is closed. As
    $$\pi \mapsto \int c \, d\pi$$
    is lower semicontinuous as noted in the proof of
    Lemma~\ref{lemma:V_is_lsc},
    $$C(\mu, \nu, \pi)= \int c \, d\pi + \iota_{G}(\mu, \nu, \pi)$$
    is lower semicontinuous, where $\iota_G=0$ on $G$ and $+\infty$ otherwise. 
\end{proof}
Now we can state and prove the first main result of this section.
\begin{theorem}\label{theo:selection_kantorovich}
     Let $(E, {\mathcal E}, m)$ be a $\sigma$-finite measure space, $X,Y$ Polish spaces, $c: X\times Y \to [0,\infty]$ a lower semicontinuous cost function and $\psi: E  \to \PP(X)\times \PP(Y)$
     a ${\mathcal E}$-measurable function.
     Then there exists a ${\mathcal E}$-measurable function $\psi': E \to \PP(X\times Y)$ so that $m$-almost everywhere $\psi'(\lambda)$ is a solution to the optimal transport problem for the marginals $\psi(\lambda)$. 
\end{theorem}

\begin{proof}
    We want to apply Corollary~\ref{cor:selection} to the set
    $$ G:= \{ (\mu, \nu, \pi) \in \PP(X)\times \PP(Y) \times \PP(X\times Y) : \pi \in \cpl(\mu, \nu)\} $$ 
    and the extended cost functional $C$. By Lemma~\ref{lemma:couplings_closed} $G$ is closed and thus Borel. By Lemma~\ref{lemma:extended_cost_is_lsc}, $C$ is lower semicontinuous.    
    Further, $\inf_{\pi \in \PP(X \times Y)}C(\mu, \nu, \pi)= V_c(\mu,\nu)$ and hence Borel measurable by Lemma~\ref{lemma:V_is_lsc}. Lastly, the existence of minimizers is guaranteed by \cite[Theorem 4.1]{Vi09}, which concludes the proof.
\end{proof}

\begin{remark}
   The conclusion of the preceding theorem can be strengthened: there exists a Borel selector of optimal transport plans on the entire space $\PP(X) \times \PP(Y)$. Indeed, the set
   $$\{(\mu, \nu, \pi): \pi \in \cpl(\mu, \nu), C(\mu, \nu, \pi)=V_c(\mu, \nu)\}$$
   is Borel and the sections at $(\mu, \nu)$ are nonempty and compact for all $(\mu, \nu)\in \PP(X)\times \PP(Y)$. Therefore, the Arsenin--Kunugui uniformization theorem \citep[Theorem~18.18]{Ke95} gives a Borel uniformization. A similar argument can be applied to the selection of Wasserstein barycenters and optimal randomized stopping times in Theorems~\ref{theo:selection_barycenters} and~\ref{theo:skorokhod}, respectively.

   We choose the Jankov--von Neumann uniformization theorem throughout as it allows a more uniform presentation: sections are not guaranteed to be countable unions of compact sets in all settings that we consider, and a selector that is optimal almost everywhere suffices for our purposes.
\end{remark}

This approach can also be used to show a selection result for optimal transport maps, given that Kantorovich-optimal transport plans of Monge type exist. We first select such transport plans, and then recover optimal maps by disintegration.

\begin{theorem}
     Let $(E, {\mathcal E}, m)$ be a $\sigma$-finite measure space, $X,Y$ Polish spaces, $c: X\times Y \to [0,\infty]$ a lower semicontinuous cost function and $\psi:  E \to \PP(X)\times \PP(Y)$ a ${\mathcal E}$-measurable function. Assume that for each $\lambda \in E$ there exists a Kantorovich-optimal transport plan of Monge type for the marginals $\psi(\lambda)$.
     Then there exists a jointly measurable function $\psi': E \times X \to  Y$ so that $m$-almost everywhere $\psi'(\lambda, \cdot)$ is a solution to the Monge optimal transport problem for the marginals $\psi(\lambda)$. 
\end{theorem}

\begin{proof}
    The extended cost function $C$ and the value function $V_c$ are Borel by Lemma~\ref{lemma:extended_cost_is_lsc} and Lemma~\ref{lemma:V_is_lsc}.    Define
    $$ G:= \{ (\mu, \nu, \pi) \in \PP(X)\times \PP(Y) \times \PP(X\times Y) : \pi \in \cpl(\mu, \nu)\}.$$
    By Lemma~\ref{lemma:couplings_closed}, $G$ is closed.
    By \cite[Proposition 5.5]{BaBeEdPi17}, the set $M \subseteq \PP(X \times Y)$ of transport plans of Monge-type is Borel.
    Consequently, the set
    $$A:= \left\{(\mu,\nu,\pi) \in G:\ C(\mu, \nu,\pi)=V_c(\mu,\nu),\ \pi \in M \right\}$$
    is Borel.
    By the Jankov–von Neumann uniformization theorem (Theorem~\ref{theo:jankov}) there exists a uniformization $\tilde f: \proj_{\PP(X) \times \PP(Y)}A \to \PP(X \times Y)$ so that $$\int c \,  d\tilde f(\mu, \nu)= V_c(\mu, \nu) $$
    and $\tilde f(\mu, \nu) \in \cpl(\mu, \nu)$ is of Monge-type for each $(\mu, \nu) \in \proj_{\PP(X) \times \PP(Y)}A$. Note that by assumption, $\psi(E) \subseteq \proj_{\PP(X) \times \PP(Y)}A$, and thus $\tilde f \circ \psi$ is well defined. Let $\sigma \in \PP(X\times Y)$ be arbitrary and define the $\A(\PP(X) \times \PP(Y))$-measurable function 
    $$ f: \PP(X) \times \PP(Y) \to \PP(X\times Y), \quad f(\mu, \nu)= \begin{cases} 
    \tilde f(\mu, \nu), &(\mu, \nu) \in \proj_{\PP(X)\times \PP(Y)}A, \\
    \sigma, &\text{otherwise.}
    \end{cases} $$
    By Lemma~\ref{lemma:appendix_2}, the function $ f \circ \psi$ is measurable with respect to the $m$-completion of ${\mathcal E}$. By Lemma~\ref{appendix:measurable_version} there exists a ${\mathcal E}$-measurable function $\hat \psi: E \to \PP(X \times Y)$, so that $\hat \psi=f\circ \psi$ $m$-almost everywhere.
     By \cite[Corollary 1.26]{Ka17} we know that there exists a measurable function ${\kappa}: \PP(X \times Y) \times X \to \PP(Y)$ so that for each $\pi \in \PP(X \times Y)$, the map ${\kappa}(\pi,\cdot)$ is a disintegration of $\pi$ with respect to its first marginal.
    Finally, note that the Dirac embedding
    $$\delta: Y \to \PP(Y), \qquad y \mapsto \delta_y$$
    is a homeomorphism onto its image. Let $y_0 \in Y$ be arbitrary and set
    \begin{equation}
    r: \PP(Y) \to Y, \qquad r(\rho):=
    \begin{cases} 
    y, &\rho=\delta_y, \\
    y_0, &\text{otherwise.}
    \end{cases} 
    \end{equation}
   Note that the Dirac measures are closed in $\PP(Y)$, and the inverse of the Dirac embedding is continuous onto its image. Together, this implies that $r$ is Borel. Given the fact that $\hat\psi(\lambda)$ is $m$-almost surely of Monge-type, and hence its disintegration is $\hat\psi(\lambda)_X$-almost everywhere a Dirac, for $m$-almost every $\lambda$, we finish the proof by setting
    \[\psi':= r \circ {\kappa} \circ (\hat \psi, \id_X).\qedhere \]
\end{proof}

If optimal transport plans exist for almost every parameter, the lower semicontinuity assumption on the cost can be relaxed to Borel measurability.

\begin{theorem}\label{theo:selection_determinacy}
  Let $(E, {\mathcal E}, m)$ be a $\sigma$-finite measure space, $X,Y$ Polish spaces, $c: X\times Y \to [0,\infty]$ a Borel measurable cost function and $\psi:  E \to  \PP(X)\times \PP(Y)$ a ${\mathcal E}$-measurable function. Assume that for $m$-almost every $\lambda \in E$, there exists a solution to the optimal transport problem given by the marginals $\psi(\lambda)$. Then there exists a ${\mathcal E}$-measurable function $\psi': E \to \PP(X\times Y)$ so that $m$-almost everywhere $\psi'(\lambda)$ is a solution to the optimal transport problem for the marginals $\psi(\lambda)$. 
\end{theorem}

\begin{proof}
For a set-valued map $F:U\to 2^V$, its graph is the set
$\{(u,v)\in U\times V:v\in F(u)\}$. The closed set $G$ from
Lemma~\ref{lemma:couplings_closed} is the graph of the coupling map
\[
    \cpl:\PP(X)\times\PP(Y)\to\mathrm{CL}\bigl(\PP(X\times Y)\bigr),
    \qquad (\mu,\nu)\mapsto\cpl(\mu,\nu),
\]
where $\mathrm{CL}(S)$ denotes the family of closed subsets of $S$.
The map $I_c:\pi\mapsto\int c\,d\pi$ is Borel. Since $G$ is closed,
the graph $A$ of the set-valued map
\[
    ((\mu,\nu),a)\mapsto
    \bigl\{\pi\in\cpl(\mu,\nu):I_c(\pi)\leq a\bigr\}
\]
is Borel in
$\bigl(\PP(X)\times\PP(Y)\times[0,\infty]\bigr)\times\PP(X\times Y)$.

By the Jankov--von Neumann uniformization theorem
\citep[Theorem~29.9]{Ke95}, there is a universally measurable
uniformization $\tilde f$ of $A$. Its domain is the projection of $A$
onto $\PP(X)\times\PP(Y)\times[0,\infty]$, which is analytic and hence
universally measurable. Extending $\tilde f$ by a fixed
$\sigma\in\PP(X\times Y)$ gives a universally measurable map
$f:\PP(X)\times\PP(Y)\times[0,\infty]\to\PP(X\times Y)$. Define
\[
    h(\mu,\nu):=f\bigl((\mu,\nu),V_c(\mu,\nu)\bigr).
\]
By Lemma~\ref{lemma:universal_value}, the map
$(\mu,\nu)\mapsto\bigl((\mu,\nu),V_c(\mu,\nu)\bigr)$ is universally
measurable. Thus $h$ is universally measurable, since this property is
preserved under composition \citep[Exercise~29.10]{Ke95}. Whenever the
infimum is attained, the section of $A$ at
$((\mu,\nu),V_c(\mu,\nu))$ is the set of optimal plans for $(\mu,\nu)$.
Therefore $h(\mu,\nu)$ is optimal at every such pair.

By assumption, $\tilde\psi':=h\circ\psi$ is optimal for the marginals
$\psi(\lambda)$ for $m$-almost every $\lambda$.
Lemma~\ref{lemma:appendix_2} shows that $\tilde\psi'$ is measurable
with respect to the $m$-completion of $\mathcal E$. 
Lemma~\ref{appendix:measurable_version} then gives the desired $\mathcal E$-measurable map
$\psi'$.
\end{proof}

\section{Wasserstein Barycenters}
For a metric space $(M,d)$, let $\PP_2(M)$ denote the set of all probability measures with finite second moment. We equip this space with the Wasserstein $2$-distance given by
$$\W_2(\mu, \nu) = \left( \inf_{\pi \in \cpl(\mu, \nu)} \int d(x,y)^2 \, \pi(dx,dy) \right)^{\frac{1}{2}},$$
for $\mu, \nu \in \PP_2(M)$.
Let $\mu \in \PP_2( \PP_2(\R^d))$ for some $d \in \N$. In analogy to the definition of barycenters of mass distributions in $\R^d$, a barycenter of $\mu$ is defined as a minimizer over $\PP_2(\R^d)$ of the functional

$$ \nu \mapsto \int \W^2_2(\sigma, \nu) \, \mu(d\sigma),$$
where $\W_2^2$ is the squared Wasserstein $2$-distance. The following theorem guarantees that under a compactness assumption, we can select barycenters in a measurable way.

\begin{theorem}\label{theo:selection_barycenters}
Let $(E,{\mathcal E}, m)$ be a $\sigma$-finite measure space, ${K}$ a nonempty, compact and convex subset of $\R^d$ and  $\psi: E \to \PP(\PP({K})) $ a ${\mathcal E}$-measurable function. Then there exists a ${\mathcal E}$-measurable function $\psi': E \to \PP({K})$ so that $m$-almost everywhere $\psi'(\lambda)$ is a barycenter of $\psi(\lambda)$.
\end{theorem}

\begin{proof}
    We first show that 
    $$\inf_{\nu \in \PP_2(\R^d)} \int \W^2_2 (\sigma, \nu) \, \mu(d\sigma)= \inf_{\nu \in \PP({K})} \int \W^2_2 (\sigma, \nu) \, \mu(d\sigma),$$
    for every $\mu \in \PP(\PP({K})).$ In particular, every minimizer over $\PP({K})$ is a barycenter.  To see this, note that the metric projection $P_{K}: \R^d \to {K}$ exists by \cite[Theorem 3.14]{BaCo11} and, as ${K}$ is nonempty, closed and convex, is $1$-Lipschitz \cite[Proposition 4.8]{BaCo11}. Now let $\sigma \in \PP({K}), \nu \in \PP_2(\R^d)$ and $\pi \in \cpl(\sigma, \nu)$. Then $(P_{K})_\#\nu \in \PP({K})$ and $(\id, P_{K})_\#\pi \in \cpl\bigl(\sigma,(P_K)_\#\nu\bigr).$ This implies that
    \begin{align*}
        \W_2^2\bigl(\sigma,(P_K)_\#\nu\bigr) &\leq \int  |x - y|^2 \, d(\id, P_{K})_\#\pi = \int \bigl|x-P_K(y)\bigr|^2 \, d\pi \\
        &= \int \bigl|P_K(x)-P_K(y)\bigr|^2 \, d\pi \leq \int |x-y|^2 \, d\pi,   
    \end{align*}
    where the second equality is due to the fact that $\sigma$ is concentrated on ${K}$ and the last inequality is due to the fact that the metric projection onto ${K}$ is $1$-Lipschitz. Now taking the infimum over all $\pi \in \cpl(\sigma, \nu)$ shows that 
    $$\W_2^2\bigl(\sigma,(P_K)_\#\nu\bigr) \leq \W_2^2(\sigma, \nu).$$ 
    This in turn implies that
    $$ \int \W_2^2\bigl(\sigma,(P_K)_\#\nu\bigr) \, \mu(d\sigma) \leq \int \W_2^2(\sigma, \nu) \, \mu(d\sigma),$$
    for every $\mu \in \PP( \PP({K})).$ Since $(P_{K})_\#\nu \in \PP({K})$ for every $\nu \in \PP_2(\R^d),$ this shows that
    $$ \inf_{\nu \in \PP({K})} \int \W^2_2 (\sigma, \nu) \, \mu(d\sigma) \leq \inf_{\nu \in \PP_2(\R^d)} \int \W^2_2 (\sigma, \nu) \, \mu(d\sigma).$$
    The other direction is trivial, as $\PP({K}) \subseteq \PP_2(\R^d)$, which proves the claim.
    
    Next, we want to apply Corollary~\ref{cor:selection} to the Borel set $\PP(\PP({K}))\times \PP({K})$ and $J$ to select barycenters. Define the map
    $$J: \PP(\PP({K})) \times \PP({K}) \to \R, \qquad J(\mu, \nu):= \int_{\PP({K})} \W^2_2(\sigma, \nu) \, \mu(d\sigma).$$
    As ${K}$ is compact, $\PP({K})$ and $\PP(\PP({K}))$ are compact by Prokhorov's theorem. Moreover, the  weak topology and the topology induced by $\W_2$ coincide on $\PP({K})$. Therefore,
    $$(\sigma, \nu) \mapsto \W^2_2(\sigma, \nu)$$
    is bounded and continuous, which implies that $J$ is continuous.
    Define the value function
    $$J': \PP(\PP({K})) \to \R, \qquad J'(\mu)=\inf_{\nu \in \PP({K})} J(\mu, \nu).$$
    Since $\PP({K})$ is compact and $J$ is continuous, the infimum is attained for every $\mu \in \PP(\PP({K}))$.

    Lastly, we show that $J'$ is continuous. Let $(\mu_n)_{n\in \N}$ be a sequence in $\PP(\PP({K}))$ so that $\lim_{n\in \N} \mu_n=\mu$ for some $\mu \in \PP(\PP({K}))$ and choose a sequence $( \nu_n)_{n\in \N}$ in $\PP({K})$ so that $\nu_n$ is a minimizer for $J(\mu_n, \cdot)$ for each $n\in \N$. Now choose a subsequence $(\mu_{n_k})_{k\in \N}$ so that 
    $\liminf_{n\in \N} J'(\mu_n) = \lim_{k \in \N} J'(\mu_{n_k}).$
    As $\PP({K})$ is compact, there exists a further subsequence of  $(\nu_{n_k})_{k\in \N}$ so that  $(\nu_{n_k} )_{k \in \N}$ converges to some $\nu \in \PP({K})$ along that subsequence. Relabeling $(\mu_{n_k})_{k \in \N}$ and $(\nu_{n_k})_{k \in \N}$ if necessary, we may assume that $\lim_{k \in \N} \nu_{n_k} =\nu$. Due to the continuity of $J$ we find that
    $$\liminf_{n \in \N} J'(\mu_n) = \lim_{k \in \N} J'(\mu_{n_k}) = \lim_{k \in \N} J(\mu_{n_k}, \nu_{n_k}) = J(\mu, \nu) \geq J'(\mu).$$
    On the other hand, if $\nu^* \in \PP({K})$ is a minimizer of $J(\mu, \cdot)$ then
    $$\limsup _{n \in \N}J'(\mu_n) \leq \lim_{n\in \N} J(\mu_n, \nu^*) = J(\mu, \nu^*)=J'(\mu).$$
    Hence, $\lim_{n \in \N} J'(\mu_n) = J'(\mu).$
    Thus all assumptions for Corollary~\ref{cor:selection} are satisfied, which concludes the proof.
\end{proof}

\section{Skorokhod Embeddings}
Let $(\Omega, \G, (\G_t)_{t \geq 0}, \P)$  be a probability space supporting a $\G_0$-measurable  random variable $U \sim \text{Unif}[0,1]$ and  a Brownian motion $B$ with respect to $(\G_t)_{t\geq 0}$. In a slight abuse of notation, we use the same notation $B$ for the associated map
$$B: \Omega \times \R_+ \to \R, \qquad B(\omega, t):=B_t(\omega). $$
We denote 
$$\PP_{2,0}(\R):=\left\{\mu\in\PP(\R):\int x\,d\mu=0,\ \int x^2\,d\mu<\infty\right\}$$
equipped with the $2$-Wasserstein topology. The Skorokhod embedding problem is to find a stopping time $\tau$ on $\Omega$ so that 
$$B_\tau \sim \mu, \qquad \E[\tau] < \infty.$$
There are many constructions of such stopping times, and several can be characterized as solutions of optimization problems. For this, define
$$S:=\{ (f, s): f:[0,s] \to \R \text{ is continuous}, f(0)=0\}$$
and equip this space with the metric
$$d_S((f,s),(g,t)):= \max \left(t-s,\ \sup_{0\leq u\leq s}\bigl|f(u)-g(u)\bigr|,\ \sup_{s\leq u\leq t}\bigl|g(u)-f(s)\bigr|\right),$$
for $(f,s),(g,t) \in S, s\leq t$. For $t < s$ the metric $d_S$ is defined symmetrically. This makes $S$ into a Polish space. Let $\gamma: S \to [0, +\infty]$ be measurable.
The optimal Skorokhod embedding problem associated to $\gamma$ consists in minimizing
$$\E\Bigl[\gamma\bigl((B_t)_{t\leq \tau},\tau\bigr)\Bigr]$$
over all solutions $\tau$ of the Skorokhod embedding problem for $\mu \in \PP_{2,0}(\R)$. We denote its value by
$$P_\gamma(\mu):= \inf \left\{  \E\Bigl[\gamma\bigl((B_{t})_{t \leq \tau},\tau\bigr)\Bigr] : \tau \text{ solves the Skorokhod embedding problem}\right\}.$$

Following the setting introduced in \cite{BeCoHu17}, we further introduce the canonical space $C_0(\R_+)$ of continuous real-valued functions on $\R_+$ starting at the value $0$, equipped with the topology of uniform convergence on compact sets. We denote the Wiener measure by $\W$. The canonical process $(\hat B_t)_{t\geq 0}$ on $C_0(\R_+)$ is given by $\hat B_t(\omega):= \omega(t)$. By a mild abuse of notation, we denote by $\hat B$ both the Brownian motion and its associated evaluation map
$$\hat B: C_0(\R_+)\times \R_+ \to \R, \qquad \hat B(\omega, t):=\hat B_t(\omega).$$
We denote by $(\F_t)_{t\geq0}$ the natural or canonical filtration and by $(\F^a_t)_{t\geq 0}$ its usual augmentation. Further, we introduce the map $r$ given by
$$r: C_0(\R_+) \times \R_+ \to S, \qquad r(\omega, t):= \bigl(\omega_{\restriction [0,t]},t\bigr).$$
It is easy to check that $r$ is continuous.

For a Polish space $X$, we denote the set of all subprobability measures on $X$ by $\PP^{\leq1}(X)$. A randomized stopping time is a measure 
$$\xi \in \PP^{\leq1}(C_0(\R_+) \times \R_+)$$
admitting a disintegration
$$\xi(d\omega,dt)=\xi_\omega(dt) \W(d\omega),$$
 where $\xi_\omega \in \PP^{\leq1}(\R_+)$ for $\W$ almost every $\omega$, so that the map
 $$\omega \mapsto \xi_\omega([0,t])$$
 is $\F_t^a$-measurable for every $t\geq 0$. The set of all randomized stopping times is denoted by $\text{RST}$. The set of all randomized stopping times with mass $1$ is denoted by $\text{RST}^1$.

For $\mu \in \PP_{2,0}(\R) $ we denote by $\text{RST}(\mu) $ the set
$$\text{RST}(\mu) := \left\{\xi\in\text{RST}^1:\hat B_\#\xi=\mu\text{ and }\int t\,\xi(d\omega,dt)<\infty\right\}.$$
As $\hat B_t^2-t$ is a martingale, this implies in particular that 
$$\int t \, \xi (d\omega, dt) = \int x^2 \, \mu(dx),$$
see \cite[Lemma 3.12]{BeCoHu17}.

We define the extended cost functional $C: \PP_{2,0}(\R) \times \text{RST}^1 \to \R \cup \{+\infty\}$ by
\begin{equation} \label{skorokhod:C}
   C(\mu, \xi)=
\begin{cases} 
\int \gamma \circ r \, d\xi, &\xi \in \text{RST}(\mu), \\
+\infty, &\text{otherwise.}
\end{cases} 
\end{equation}
With this, we have by \cite[Lemma 3.11]{BeCoHu17} that 
\begin{equation}\label{skorokhod:P_gamma}
     P_\gamma(\mu)= \inf_{\xi \in \text{RST}^1} C(\mu,\xi).
\end{equation}

\begin{lemma}\label{lemma:convergence_for_skorokhod}
    Let $(\mu_n)_{n\in \N}$ be a sequence of measures in $\PP_{2,0}(\R)$. Assume that $\lim_{n \to \infty} \mu_n =\mu$ in the $2$-Wasserstein sense for a measure $\mu \in \PP_{2,0}(\R)$. Let further $(\xi_n)_{n \in \N}$ be a sequence of randomized stopping times so that $\xi_n \in \text{RST}(\mu_n)$ for all $n \in \N$. Then there exists a subsequence $(\tilde \xi_n)_{n \in \N}$ of $(\xi_n)_{n\in \N}$ that converges weakly to a randomized stopping time $\xi \in \text{RST}(\mu).$
\end{lemma}

\begin{proof}
    We first show that $(\xi_n)_{n\in \N}$ is tight by showing that the sets of marginals of $(\xi_n)_{n \in \N}$ are tight. As $(\xi_n)_{n \in \N}$ are randomized stopping times, they all have the same marginal ${\W}$ as their marginal on $C_0(\R_+)$. As for the other marginal, denote by $T: C_0(\R_+)\times \R_+ \to \R_+$ the projection given by $T(\omega, t):=t$. As $\xi_n \in \text{RST}(\mu_n)$ for all $n \in \N$ and $\lim_{n \to \infty} \mu_n =\mu$ in the $2$-Wasserstein sense, we have
    $$ \lim_{n \to \infty}\E_{\xi_n}[T]=\lim_{n\to \infty} \int x^2 \, d\mu_n = \int x^2 \, d\mu ,$$
    which implies that $\sup \E_{\xi_n}T < \infty$ and hence the marginals of $(\xi_n)_{n \in \N}$ on $\R_+$ are also tight. Thus,  $(\xi_n)_{n \in \N}$ admits a subsequence $(\tilde \xi_n)_{n \in \N}$ that converges weakly to measure $\xi \in \PP( C_0(\R_+) \times \R_+)$. As the set of randomized stopping times is closed (\cite[Corollary 3.10]{BeCoHu17}), $\xi$ is a randomized stopping time. Next we show that $\hat B_\#\xi=\mu$. As the evaluation map $\hat B(\omega, t)=\omega(t)$ is continuous,
    $$\lim_{n \to \infty} \hat B_\# \tilde \xi_n = \hat B_\# \xi,$$
    weakly. But at the same time
    $$\lim_{n \to \infty} \hat B_\#\xi_n =\lim_{n \to \infty} \mu_n =\mu$$
    in the $2$-Wasserstein sense, and thus in particular also
    $$\lim_{n \to \infty} \hat B_\#\tilde\xi_n =\mu.$$
     By uniqueness of limits, we get
    $$ \hat B_\#\xi = \mu.$$
It remains to show that $\xi \in \text{RST}(\mu)$. By the portmanteau theorem,
$$ \E_\xi [T] \leq \liminf_{n \to \infty} \E_{\tilde \xi_n} [T] < \infty.$$
By \cite[Lemma 3.12]{BeCoHu17} this implies that $\E_\xi [T] = \int x^2 \, d\mu$, and hence $\xi \in \text{RST}(\mu)$.
\end{proof}

\begin{theorem}\label{theo:skorokhod}
    Let $(E,{\mathcal E}, m)$ be a $\sigma$-finite measure space and $\psi: E \to \PP_{2,0}(\R) $ a ${\mathcal E}$-measurable function.
    Let further $\gamma: S \to [0,+\infty]$ be lower semicontinuous.
    Then there exists a ${\mathcal E}$-measurable function $\psi': E \to \text{RST}^1$ so that $m$-almost everywhere $\psi'(\lambda)$ is a minimizer of $P_\gamma(\psi(\lambda)).$
\end{theorem}

\begin{proof}
    We want to apply Corollary~\ref{cor:selection} to the set
    $G:= \{(\mu, \xi) \in \PP_{2,0}(\R)\times \text{RST}^1: \xi \in \text{RST}(\mu)\}$
    and the extended cost functional $C$ defined in \eqref{skorokhod:C}. We first note that $G$ is closed, and thus Borel. Indeed, let $(\mu_n, \xi_n)_{n\in \N}$ be such, that $(\mu_n, \xi_n) \in G$  for all $n \in \N$ and $\lim_{n\in \N}(\mu_n, \xi_n)=(\mu, \xi)\in \PP_{2,0}(\R)\times \text{RST}^1$. By Lemma~\ref{lemma:convergence_for_skorokhod} there exists a subsequence of $(\xi_n)_{n\in \N}$ that converges to a randomized stopping time $\xi' \in \text{RST}(\mu)$. By uniqueness of limits, $\xi = \xi'$, and hence $\xi \in \text{RST}(\mu)$.   Further note that $\text{RST}^1$ equipped with the weak topology is a Polish space, as it is a closed subset of $\PP(C_0(\R_+) \times \R_+)$ by \cite[Corollary 3.10]{BeCoHu17}.
    
    Next, we show that $C$ is measurable by establishing that it is lower semicontinuous. To this end, let $(\mu_n,\xi_n)_{n \in \N}$ be a sequence of elements in $\PP_{2,0}(\R) \times \text{RST}^1$ converging to a pair $(\mu,\xi) \in \PP_{2,0}(\R) \times \text{RST}^1$. Without loss of generality we may assume that $\liminf_{n \to \infty}C(\mu_n, \xi_n) < \infty$ and $\xi_n \in \text{RST}(\mu_n)$ for all $n \in \N$. By Lemma~\ref{lemma:convergence_for_skorokhod} we know that $\xi \in \text{RST}(\mu)$. Note that for a suitable subsequence $(\mu_{n_k}, \xi_{n_k})_{k \in \N}$ we have
    $$\lim_{k \to \infty} C(\mu_{n_k}, \xi_{n_k}) =\liminf_{n \to \infty} C(\mu_n, \xi_n). $$
    Relabeling the sequence $(\mu_n, \xi_n)_{n \in \N}$ if necessary, we may assume the same for this sequence. 
    As $r$ is continuous and  $\gamma$ is non-negative and lower semicontinuous, we further know by the Portmanteau theorem that $\xi \mapsto \int \gamma \circ r\, d\xi$ is lower semicontinuous, and thus
    $$C(\mu, \xi)= \int \gamma \circ r\, d\xi \leq \liminf_{n \to \infty} \int \gamma \circ r\, d\xi_n = \liminf_{n \to \infty} C(\mu_n, \xi_n).$$
    By \cite[Theorem 3.14]{BeCoHu17}, the set $\text{RST}(\mu)$ is non-empty and compact for every $\mu \in \PP_{2,0}(\R)$.
    By lower semicontinuity, the functional
    $$\xi \mapsto \int \gamma \circ r \, d\xi$$
    attains its infimum on $\text{RST}(\mu)$, proving the existence of optimal randomized stopping times for each $\mu \in \PP_{2,0}(\R).$
   It remains to show that $ P_\gamma(\cdot)$ is measurable, which we once more do by establishing the stronger property of being lower semicontinuous. Hence, let $(\mu_n)_{n \in \N}$ be a sequence in $\PP_{2,0}(\R)$ converging in the $2$-Wasserstein metric to a measure $\mu \in \PP_{2,0}(\R)$. We know that for a suitable subsequence  $(\mu_{n_k})_{k \in \N}$ we have
    $$\lim_{k \to \infty} P_\gamma(\mu_{n_k}) = \liminf_{n \to \infty} P_\gamma(\mu_n).$$
    Relabeling the sequence $(\mu_n)_{n \in \N}$ if necessary, we may assume the same holds for this sequence.
    Choose $(\xi_n)_{n\in \N}$ such that $ \xi_n \in \text{RST}(\mu_n)$ and 
    $$\int \gamma \circ r \, d\xi_n = P_\gamma(\mu_n),$$
    for each $n \in \N$. By Lemma~\ref{lemma:convergence_for_skorokhod} there exists a subsequence of $(\xi_n)_{n\in \N}$ that converges weakly to a randomized stopping time $\xi \in \text{RST}(\mu)$. Using once again that $\xi \mapsto \int \gamma \circ r\, d\xi$ is lower semicontinuous, we calculate
    $$P_\gamma(\mu) \leq \int \gamma \circ r \, d\xi \leq \liminf_{n \to \infty} \int \gamma \circ r \, d\xi_n = \liminf_{n \to \infty} P_\gamma(\mu_n).$$
  Thus, all assumptions for the measurable selection corollary are satisfied, yielding the desired measurable selector.
\end{proof}
Next, we want to use the preceding theorem to show a measurable selection result for ordinary stopping times. We first establish a measurability property of parametrized random variables.

\begin{lemma}\label{lemma:jointly_measurable_to_L0}
    Let $(S, \F)$ be a measurable space, $\Omega$ a Polish space, $\mu \in \PP(\Omega)$ and 
    $$f: S \times \Omega \to \R$$
    jointly measurable. Then the map $\phi$ given by
    $$\phi: S \to L^0(\Omega, \mu), \qquad \phi(s):= [f(s, \cdot)]$$
    is measurable.
\end{lemma}

\begin{proof}
Take a bounded metric $d_\R$ on $\R$ which is compatible with the topology of $\R$. Then, the metric $d_{L^0}$ given by
$$d_{L^0}: L^0(\Omega, \mu) \times L^0(\Omega, \mu) \to \R, \qquad d_{L^0}(X,Y):= \int_{\Omega} d_\R(X,Y) \, d\mu$$
metrizes convergence in probability.  As $f$ is jointly measurable, the map
$$(s, \omega) \mapsto d_\R(f(s, \omega), g(\omega))$$
is $(\F \otimes \B(\Omega))$-measurable and non-negative for every measurable $g: \Omega \to \R$. The standard measurability theorem for parameter-dependent integrals then implies that also the map $\psi_g$ defined by
$$\psi_g: S \to \R, \qquad \psi_g(s):=\int_\Omega d_\R(f(s,\omega), g(\omega)) \, d\mu = d_{L^0}\bigl([f(s,\cdot)],[g]\bigr)$$
is measurable. Note that if $[g]=[g']$ for measurable $g,g': \Omega \to \R$, then $\psi_g=\psi_{g'}$. In a slight abuse of notation, we hence also write $\psi_Z$ for $Z\in L^0(\Omega, \mu)$.
By \cite[Proposition 19.4-19.5]{Ke10} the space $L^0(\Omega, \mu)$ is separable. Let $(Z_i)_{i \in \N}$ be a countable, dense collection of elements of $L^0(\Omega, \mu)$. Then the open balls $(B_q(Z_i))_{i \in \N, q \in \Q_{>0}}$ form a countable basis for $L^0(\Omega, \mu)$ and hence generate the Borel sets on $L^0(\Omega, \mu)$. But
$$\phi^{-1}(B_q(Z_i)) = \{s : d_{L^0}\bigl([f(s,\cdot)],Z_i\bigr) < q\} = \psi^{-1}_{Z_i}\bigl([0,q)\bigr),$$
which belongs to $\F$.
\end{proof}
In the following theorem, we establish a selection result in which we want to measurably select stopping times. As the set of stopping times does not possess a natural Borel $\sigma$-algebra, we identify almost surely finite stopping times, up to almost sure equality, with elements of $L^0$, allowing us to make precise statements about measurability. 
\begin{theorem}\label{theo:selection_stopping_times}
    Let $(\Omega, \G, (\G_t)_{t \geq 0}, \P)$ be a stochastic basis supporting a Brownian motion $B$ with respect to $(\G_t)_{t\geq0}$ and a $\G_0$-measurable random variable $U \sim \text{Unif}[0,1]$. Let $(E, {\mathcal E}, m)$ be a $\sigma$-finite measure space and let 
    $$\psi:  E \to \PP_{2,0}(\R)$$
    be ${\mathcal E}$-measurable. Let further $\gamma: S \to [0, +\infty]$ be lower semicontinuous. Then there exists a ${\mathcal E}$-measurable function 
    $$\psi': E \to L^0(\Omega, \P )$$
    so that, for $m$-almost every $e \in E$, the equivalence class $\psi'(e)$ admits a representative $\tau_e$ which is an almost surely finite stopping time satisfying
    $$B_{\tau_e} \sim \psi(e), \qquad \E_\P[\tau_e] < \infty$$
    and
    $$\E \Bigl[\gamma\bigl((B_t)_{t \leq \tau_e},\tau_e\bigr)\Bigr]= P_\gamma (\psi(e)).$$ 
\end{theorem}

\begin{proof}
By Theorem~\ref{theo:skorokhod} there exists a ${\mathcal E}$-measurable function $\hat \psi: E \to   \text{RST}^1$ so that $m$-almost everywhere $\hat \psi(e)$ is a minimizer of $P_\gamma(\psi(e))$. We separate this proof into two steps: we first construct a measurable function $\phi: \text{RST}^1 \to L^0(C_0(\R_+) \times [0,1], \W \otimes \lambda) $, where $\W$ is the Wiener measure and $\lambda$ is the Lebesgue measure on $[0,1]$, so that for every $\xi \in \text{RST}^1$, the image $\phi(\xi)$ has a representative $\rho_\xi$ that is a stopping time so that
$$\hat B_{\rho_\xi} \sim  \hat B_\#\xi$$
and
$$ \int_{C_0(\R_+) \times \R_+} \gamma \circ r \, d\xi = \E_{\W \otimes \lambda}\Bigl[\gamma\bigl((\hat B_t)_{t\leq \rho_\xi},\rho_\xi\bigr)\Bigr].$$
In a second step, we construct a measurable function $\phi': L^0(C_0(\R_+) \times [0,1], \W \otimes \lambda) \to L^0(\Omega, \P) $ with the following property: If  $f \in L^0(C_0(\R_+) \times [0,1], \W \otimes \lambda)$ has a representative $\rho$ that is stopping time on $C_0(\R_+) \times [0,1]$, then $\phi'(f)$ has a representative $\tau_\rho$ so that 
$$\hat B_\rho \sim B_{\tau_\rho}$$
and 
$$\E_{\W \otimes \lambda}\Bigl[\gamma\bigl((\hat B_t)_{t\leq \rho},\rho\bigr)\Bigr] = \E_\P\Bigl[\gamma\bigl((B_t)_{t\leq \tau_\rho},\tau_\rho\bigr)\Bigr].$$

For the first step, begin by setting 
$$\rho_\xi(\omega, u):=\inf \{t \geq 0: \xi_\omega([0,t]) \geq u\}$$
for $\xi \in \text{RST}^1$, where $\xi_\omega$ is a disintegration of $\xi$ with respect to $\omega$. As $\rho_\xi$ is almost surely finite, we interpret $[\rho_\xi]$ as the equivalence class of any measurable, finite-valued modification of $\rho_\xi$. Then set
$$\phi(\xi):=[\rho_\xi].$$
  By \cite[Theorem 3.8]{BeCoHu17}, $\rho_\xi$ is a stopping time for the usual augmentation of the filtration $\F_t \otimes \B([0,1])_{t\geq 0}$, where $\B([0,1])$ are the Borel sets on $[0,1]$. Since the stopping time obtained from \cite[Theorem 3.8]{BeCoHu17} is predictable with respect to the augmented filtration, \cite[Theorem~3.1]{BeCoHu17} establishes that there exists an element in the equivalence class of $\rho_\xi$ that is a stopping time with respect to the filtration $\bigl(\F_t\otimes\B([0,1])\bigr)_{t \geq 0}$. By \cite[Equation (3.9)]{BeCoHu17} we know that $\hat B_{\rho_\xi} \sim \hat B_\#\xi $ and by \cite[Lemma 3.11]{BeCoHu17} that
$$\int_{C_0(\R_+)\times \R_+} \gamma \circ r \, d\xi =\E_{\W \times \lambda}\Bigl[\gamma\bigl((\hat B_t)_{t \leq \rho_\xi},\rho_\xi\bigr)\Bigr].$$
Note, that the integral on the right is independent of particular member of the equivalence class $[\rho_\xi]$. We now show that $\phi$ is Borel measurable. By \cite[Corollary 1.26]{Ka17}  there exists a jointly Borel measurable function 
$${\kappa}: \PP(C_0(\R_+) \times \R_+) \times C_0(\R_+) \to \PP(\R_+),$$ 
so that for every $\xi \in \text{RST}^1$, the map $\omega \mapsto {\kappa}(\xi, \omega)$ is a disintegration of $\xi$ with respect to the Wiener measure. Applying \cite[Proposition 7.29]{BeSh78} to the kernel $q$ from $\PP(\R_+)\times \R_+$ to $\R_+$ given by
$$q(dx|\mu, t):= \mu(dx),$$
and the function
$$f: (\PP(\R_+)\times \R_+) \times \R_+ \to \R_+, \qquad  f((\mu, t), x):= \ind{\{x \leq t\}}$$
shows that the map
$$g:\PP(\R_+) \times \R_+ \to \R_+, \qquad g(\mu, t):= \int f((\mu, t),x) \, q(dx|\mu, t)= \mu([0,t])$$
is Borel measurable.
This implies that the map $F: \PP(C_0(\R_+)\times \R_+) \times C_0(\R_+) \times \R_+ \to [0,1] $ given by
$$ F(\xi, \omega, t):= {\kappa}(\xi, \omega)([0,t])= g\bigl(\kappa(\xi,\omega),t\bigr)$$
is also jointly Borel measurable. Consequently, the quantile
$$ R: \PP(C_0(\R_+) \times \R_+) \times C_0(\R_+) \times [0,1] \to [0,+\infty], $$
given by
$$ R(\xi, \omega, u):= \inf \{ t\geq 0: F(\xi, \omega, t) \geq u\}$$
is jointly Borel measurable as
$$ \{ R < a\} = \bigcup_{\substack{q \in \Q_+  \\ q < a}}\{(\xi, \omega, u): F(\xi, \omega, q)\geq u\},$$
for $a \in \R$.
As ${\kappa}(\xi, \cdot)$ and $\xi_\omega$ are both disintegrations of $\xi$, they agree $\W$-almost everywhere. Consequently, for every $\xi \in \text{RST}^1$,
$$\bigl[R(\xi,\cdot,\cdot)\bigr]=\bigl[\rho_\xi(\cdot,\cdot)\bigr]=\phi(\xi).$$
Moreover,   $R(\xi, \cdot, \cdot)$ is almost surely finite valued, for every  $\xi \in \text{RST}^1$. We may replace $R$ with a real-valued modification $\tilde R$, by changing all infinite values to $0$, which neither changes its equivalence class in $L^0(C_0(\R_+)\times [0,1], \W \otimes \lambda)$, nor measurability. Now, applying Lemma~\ref{lemma:jointly_measurable_to_L0} to $\tilde R$ shows that $\phi$ is Borel measurable.

For the second step, define the function $\Phi$ as 
$$\Phi: \Omega \to C_0(\R_+)\times [0,1], \qquad \Phi(\omega):= \bigl((B_t(\omega))_{t\geq 0},U(\omega)\bigr),$$
and set $\phi'$ as
$$\phi'(f):= f \circ \Phi.$$
If $f \in L^0(C_0(\R_+) \times [0,1], \W \otimes \lambda)$ has a representative $\rho$ that is a stopping time, we define the random time $\tau_\rho:= \rho \circ \Phi$.
As $U$ is $\G_0$-measurable and $B$ is a Brownian motion with respect to $(\G_t)_{t\geq 0}$, $U$ is independent of $B$ and hence 
$$  \Phi_\#\P = \W \otimes \lambda.$$
This implies that the definition of $\phi'(f)$ does not depend on the representative chosen for $f$, and thus $\phi'$ is well defined. It also follows that
$$\hat B_\rho \sim B_{\tau_\rho}$$
and
$$\E_{\W \otimes \lambda}\Bigl[\gamma\bigl((\hat B_t)_{t\leq \rho},\rho\bigr)\Bigr] = \E_\P\Bigl[\gamma\bigl((B_t)_{t\leq \tau_\rho},\tau_\rho\bigr)\Bigr]$$
whenever $\rho$ is a stopping time.
Next we show that $\tau_\rho$ is a stopping time whenever $\rho$ is a stopping time. As $\Phi$ is measurable as a map from $(\Omega, \G_t)$ to $(C_0(\R_+)\times [0,1], \F_t \times \B([0,1]))$ for every $t\geq 0$, we have that
$$\{\tau_\rho \leq t\} = \Phi^{-1}\Bigl(\rho^{-1}\bigl([0,t]\bigr)\Bigr) \in \G_t $$
for every $t\geq 0$, which shows that $\tau_\rho$ is a stopping time.
Lastly, we show that $\phi'$ is continuous, and thus also measurable. Let $f_n,f \in L^0(C_0(\R_+) \times [0,1])$ for $n \in \N$ and $\epsilon >0$, then
\begin{align*}
&\P( |\phi'(f) - \phi'(f_n)| \geq \epsilon) = \P (|f \circ \Phi - f_n\circ \Phi | \geq \epsilon ) \\
= &\P\Bigl(\Phi^{-1}\bigl(|f-f_n|\geq\epsilon\bigr)\Bigr) = \W\otimes \lambda (|f-f_n| \geq \epsilon). 
\end{align*}
Thus, $\lim_{n \to \infty} \phi'(f_n)=\phi'(f)$ in probability whenever $\lim_{n \to \infty} f_n=f$ in probability.
Now, set $\psi':= \phi' \circ \phi \circ \hat \psi$.

What is left to show, is that selected stopping times have finite expectation. Fix $e\in E$ outside the exceptional  $m$-null set, set $\mu:= \psi(e)$, and let $\xi:=\hat \psi(e) \in \text{RST}(\mu)$. Let $\tau_e$ and  $\rho_\xi$ be the corresponding stopping times. Then
$\E[\tau_e] < \infty $.
Indeed,
$$\E_\P[\tau_e] = \E_{\W \otimes \lambda}[\rho_\xi]= \int t \, \xi(d\omega, dt) = \int_\R x^2 \, \mu(dx) < \infty, $$
where the first equality follows from $\Phi_\#\P = \W \otimes \lambda$, the second follows from the quantile representation of $\xi$, and the third follows from \cite[Lemma~3.12]{BeCoHu17}.
\end{proof}

\section{Appendix}

\subsection{Analytic sets and uniformization}

We collect some fundamental and well-known properties of analytic sets here. For a more thorough discussion and proofs, see for example \cite{Ke95}.

\begin{definition}
    Let $X$ be a Polish space and $A\subseteq X$. $A$ is called analytic, if there exists a Polish space $Y$ and a continuous function $f: Y \to X$ so that $f(Y)=A$.
\end{definition}
\begin{remark}
    If $A \neq \emptyset$, then the space $Y$ can be chosen to be $ \N^\N$.
\end{remark}

There are several alternative characterizations of analytic sets, some of which are collected in the next theorem.

\begin{theorem}\label{app:equivalent_formulations_of_analytic}
    Let $X$ be a Polish space and $A\subseteq X$, then the following are equivalent:
    \begin{enumerate}
        \item $A$ is analytic.
        \item For any uncountable Polish space $Y$, there exists a Borel set $B \subseteq X \times Y$ so that $A = \proj_X(B) $.
        \item There exists a closed set $F \subseteq X \times \N^\N$ such that $A=\proj_X(F)$.
    \end{enumerate}
\end{theorem}

One important property of analytic sets is, that Borel images of analytic sets are analytic.

\begin{prop}\label{prop:analytic_images}
    Let $X,Y$ be Polish spaces, $f: X\to Y $ Borel measurable and $A \subseteq X$ analytic. Then $f(A)$ is analytic.
\end{prop}

\begin{definition}
    Let $X$ be a Polish space and $A \subseteq X$. Then $A$ is co-analytic if $A^\complement$ is analytic.
\end{definition}

There is a close relationship between analytic sets and Borel sets.

\begin{theorem}
    Let $X$ be a Polish space and $A\subset X$. $A$ is Borel if and only if it is analytic and co-analytic.
\end{theorem}

\begin{theorem}
    Let $X,Y$ be Polish spaces and $f: X \to Y$ a function. Then the following are equivalent:
    \begin{enumerate}
        \item f is Borel measurable.
        \item $\textit{graph}(f)$ is Borel.
        \item $\textit{graph}(f)$ is analytic.
    \end{enumerate}
\end{theorem}

A consequence of the Choquet capacitability theorem is the following theorem.

\begin{theorem}
    Let $X$ be a Polish space, $A \subset X$ and $\mu$ a finite Borel measure on $X$. If $A$ is analytic, then $A$ is measurable with respect to the completion of $\mu$.
\end{theorem}

The following notion of a uniformization and the Jankov–von Neumann uniformization theorem underpin a lot of results of this paper.

\begin{definition}
    Let $X,Y$ be arbitrary sets and $P \subseteq X \times Y$, then a uniformization of $P$ is a function
    $$f: \proj_X(P) \to Y$$
    such that
    $$(x, f(x)) \in P$$
    for all $x \in \proj_X(P)$.
\end{definition}

A subset of $X$ is called universally measurable if it is measurable with respect to the completion of every finite Borel measure on $X$.
We denote by
$$  \A(X):=\sigma(\{A\subseteq X:A\text{ is analytic}\}) $$
the $\sigma$-algebra generated by the analytic subsets of $X$. As analytic sets are universally measurable and the universally measurable sets form a $\sigma$-algebra, every set in $\A(X)$ is universally measurable. We further denote the $\sigma$-algebra of all Borel sets of the Polish space $X$ by $\B(X)$.

\begin{theorem}[Jankov–von Neumann uniformization theorem]\label{theo:jankov}
    Let $X,Y$ be Polish spaces and $P \subseteq X \times Y$ an analytic set. Then there exists a uniformization 
    $$h: \proj_X(P)\to Y$$
    of $P$ that is measurable between the measurable spaces $(\proj_X(P), \A(X)|_{\proj_X(P)})$ and $(Y,\B(Y))$.
\end{theorem}

For our applications, we select minimizers of several minimization problems. The following lemmas are consequences of the Jankov-von Neumann uniformization theorem applied to families of such optimization problems parametrized by a $\sigma$-finite measure space. It is useful to allow for a separate set $A$ of admissible elements, as achieving the optimal value alone is not enough to ensure admissibility, when the optimal value is $+\infty$.

\begin{lemma}\label{lemma:appendix_1}
    Let $X,Y$ be Polish spaces, $A \subseteq X \times Y$ Borel and  $f: X\times Y \to \R \cup \{+\infty\}$ be Borel measurable. Let the value function 
    $$v: X \to \R \cup \{ +\infty\}, \qquad v(x):=\inf_{y\in Y} f(x,y)$$
    be Borel measurable. If for every $x \in X$, there exists $y \in Y$, so that
    $$f(x,y)=v(x), \qquad (x,y)\in A,$$
    then the set
    $$P:=\{(x,y)\in A: f(x,y)=v(x)\}$$
    is Borel. In particular, there exists an $\A(X)$-measurable function $h: X \to Y$ so that
    $$f(x, h(x))=v(x), \qquad (x, h(x)) \in A$$
    for all $x\in X$.
\end{lemma}

\begin{proof}
     Let us introduce the auxiliary function 
    $$g: X \times Y \to (\R \cup \{+\infty\}) \times (\R \cup \{+\infty\}), \qquad g(x,y):=\bigl(f(x,y),v(x)\bigr).$$ 
    Note that  $g$ is measurable due to the assumed measurability of $f$ and $v$. Now we can express $P$ as  
    $$P=  g^{-1}\bigl(\Delta_{\R\cup\{+\infty\}}\bigr) \cap A $$
    where $\Delta_{\R \cup \{+\infty\}}$ is the diagonal of $\R \cup \{+\infty\}$. As diagonals are closed, $P$ is Borel. By Theorem~\ref{theo:jankov} we get a uniformization $h$ of $P$ which is $\A(X)|_{\proj_X(P)}$-measurable. Since we assumed that for every $x \in X$ there exists at least one $y\in Y$ so that $f(x,y)=v(x)$ and $(x,y)\in A$, we know $\proj_X(P)=X$, thus $h$ is a total function from $X$ to $Y$ and $\A(X)$-measurable.
\end{proof}

\begin{lemma}\label{lemma:appendix_2}
    Let $(E, {\mathcal E}, m)$ be a $\sigma$-finite measure space, $X,Y$ Polish spaces, functions $\psi: E \to X$ ${\mathcal E}$-measurable and  $h: X\to Y$ universally measurable. Then $h \circ \psi$ is measurable with respect to the $m$-completion of ${\mathcal E}$.
\end{lemma}

\begin{proof}
     As $m$ is $\sigma$-finite, we can find a collection of sets $(E_n)_{n\in \N}$  so that for all $n \in \N$ we have
    $$E_n \in {\mathcal E}, \qquad E= \bigcup_{n\in \N} E_n , \qquad m(E_n)< \infty. $$
    Now, for each $n \in \N$ define a finite Borel measure on $X$ by
    $$\nu_n:=\psi_{\#}(m|_{E_n}).$$
    Consider $B\in \B(Y)$ and set
    $$C:=h^{-1}(B).$$
    As $h$ is universally measurable, $C$ is universally measurable. In particular, it is measurable with respect to the $\nu_n$-completion of $\B(X)$ for each $n \in \N$. Therefore, for each $n \in \N$ there exist Borel $\tilde C_n, N_n \subseteq X$ such that
    $$C \triangle \tilde C_n \subseteq N_n, \qquad \nu_n(N_n)=0.$$
   Consequently,
   $$\bigl(\psi^{-1}(C)\cap E_n\bigr) \triangle \bigl(\psi^{-1}(\tilde C_n)\cap E_n\bigr) \subseteq \psi^{-1}(N_n) \cap E_n.$$
    Since 
    $$m\bigl(\psi^{-1}(N_n)\cap E_n\bigr) = \nu_n (N_n)=0$$
    and
    $$\psi^{-1}(\tilde C_n) \cap E_n \in {\mathcal E},$$
    it follows that $\psi^{-1}(C) \cap E_n$ is measurable with respect to the $m$-completion of ${\mathcal E}$.
    Finally,
    $$ \psi^{-1}(C) = \bigcup_{n\in \N}\bigl(\psi^{-1}(C)\cap E_n\bigr), $$
    so $\psi^{-1}(C)$ is measurable with respect to the $m$-completion of ${\mathcal E}$. As $B \in \B(Y)$ was arbitrary, this finishes the proof.     
\end{proof}

\begin{lemma}\label{appendix:measurable_version}
    Let $(E, {\mathcal E}, m)$ be a measure space, $X$  a Polish space and 
    $$f: E \to X$$
    be measurable with respect to the $m$-completion of $\mathcal E$. Then there exists a $\mathcal E$-measurable function 
    $$g: E \to X$$
    so that $m$-almost everywhere
    $$f=g.$$
\end{lemma}

\begin{proof}
    As $X$ is a Polish space, its Borel $\sigma$-algebra is countably generated. Let $(B_n)_{n \in \N}$ be such, that $B_n \in \B(X)$ and $\B(X)=\sigma(B_1,B_2, \dots)$. Since $f$ is measurable with respect to the $m$-completion of $\mathcal E$, there exist sets $A_n,N_n \in \mathcal E$ such that $m(N_n)=0$ and
    $$f^{-1}(B_n) \triangle A_n \subseteq N_n,$$
    for each $n \in \N$. Now define
    $N: = \bigcup_{n \in \N} N_n.$
    Then $N \in \mathcal E$ and $m(N)=0$. Fix $x_0 \in X$ and set
    $$g(s):=
    \begin{cases} 
    f(s), &s \in E\setminus N, \\
    x_0, &s \in N.
    \end{cases} $$
    Clearly, $m$-almost everywhere, we have $f=g$. Moreover, for each $n \in \N$, we find that
    $$g^{-1}(B_n) = 
    \begin{cases} 
    A_n \setminus N & x_0 \notin B_n, \\
    (A_n \setminus N) \cup N, & x_0 \in B_n.
    \end{cases} $$
    Thus, $g^{-1}(B_n) \in \mathcal E$ for each $n \in \N$. As the sets $(B_n)_{n \in \N}$ generate $\B(X)$, it follows that $g$ is $\mathcal E$ measurable.
\end{proof}

\begin{corollary}\label{cor:selection}
    Let $(E, {\mathcal E}, m)$ be a $\sigma$-finite measure space, $X,Y$ Polish spaces, $A \subseteq X \times Y$ Borel, functions $\psi: E \to X$ ${\mathcal E}$-measurable and  $f: X\times Y \to \R \cup \{+\infty\}$ Borel measurable. If the value function 
    $$v: X \to \R \cup \{ +\infty\}, \qquad v(x):=\inf_{y\in Y} f(x,y)$$
    is Borel measurable, and for every $x \in X$, there exists $y \in Y$ such that
    $$f(x,y) =v(x), \qquad (x,y) \in A,$$
    then there exists a ${\mathcal E}$-measurable function 
    $$\psi': E \to Y$$
    so that $m$-almost everywhere we have 
    $$f(\psi(\lambda),\psi'(\lambda))=v(\psi(\lambda)), \qquad (\psi(\lambda), \psi'(\lambda))\in A.$$
\end{corollary}

\begin{proof}
   By Lemma~\ref{lemma:appendix_1} we see that there exists a $\A(X)$-measurable function $h: X\to Y$ so that
   $$f(x,h(x))=v(x), \qquad (x, h(x))\in A,$$
   for all $x \in X$. As any $\A(X)$-measurable function is in particular universally measurable, we know by Lemma~\ref{lemma:appendix_2} that $h \circ \psi$ is measurable with respect to the $m$-completion of ${\mathcal E}$.
     By Lemma~\ref{appendix:measurable_version}, there exists a ${\mathcal E}$-measurable function
    $$\psi': E \to Y$$
    so that 
    $$ \psi'= h \circ \psi,$$
    $m$-almost everywhere. As,
    $$f(x,h(x))=v(x), \qquad (x, h(x))\in A$$
    for all $x \in X$, we get
    $$f(\psi(\lambda),\psi'(\lambda))=v(\psi(\lambda)), \qquad (\psi(\lambda),\psi'(\lambda)) \in A,$$
    $m$-almost surely.
    
\end{proof}

{\bf Acknowledgements.}
The author is grateful to Mathias Beiglböck and Manu Eder for helpful discussions throughout this work.  The author used a large language model for assistance
with exposition and technical checks. All mathematical arguments are due to the author.

This research was funded in whole or in part by the Austrian Science
Fund (FWF)  10.55776/P34743 and 10.55776/P36835 as well as the Austrian National Bank
[Jubiläumsfonds, project 18983]. For open access purposes, the author has
applied a CC BY public copyright license to any author accepted
manuscript version arising from this submission.

\end{document}